\documentclass{amsart}
\title{The Rank Stable Topology of Instantons on $\cpbar$}
\author{Jim Bryan\\
Marc Sanders}

\usepackage{amsmath,amsthm,amsfonts}

\chardef\bslash=`\\ 

\newtheorem{thm}{Theorem}[section]
\newtheorem{theorem}{Theorem}[section]

\newtheorem{conj}[thm]{Conjecture}
\newtheorem{lem}[thm]{Lemma}

\theoremstyle{definition}

\theoremstyle{remark}
\newtheorem{rem}{Remark}[section]

\newcommand{\cnums} {{\mathbf C}}          
\newcommand{\cpbar} {\overline{\mathbf {CP}}^2}
\newcommand{\dbar}{\bar\partial}

\newcommand{\notexists}{\exists\hspace{-6pt}\raisebox{1.5pt}{/}}

\newcommand{\til}[1]{\widetilde{#1}}

\newcommand{\1}{{{\mathchoice {\rm 1\mskip-4mu l} {\rm 1\mskip-4mu l}
{\rm 1\mskip-4.5mu l} {\rm 1\mskip-5mu l}}}}

\renewcommand{\P}{{\mathbf{P}}}
\newcommand{\M}{{\mathcal{M}}}
\newcommand{\E}{{\mathcal{E}}}
\renewcommand{\O}{{\mathcal{O}}}
\newcommand{\ms}[2]{\M_{#1}^{#2}}
\newcommand{\msalg}[2]{\M_{\text{alg},#1}^{#2} }
\newcommand{\linf}{\ell_{\infty }}

\begin{document}
\maketitle

\begin{abstract}
Let $\M _{k}^{n}$ be the moduli space of based (anti-self-dual) instantons
on $\cpbar $ of charge $k$ and rank $n$.  There is a natural inclusion $\M
_{k}^{n}\hookrightarrow \M_{k}^{n+1}$.  We show that the direct limit space
$\M_k^\infty$ is homotopy equivalent to $BU(k)\times BU(k)$.
Let $\linf$ be a line in the
complex projective plane and let $\til{\cnums \P}^{2}$ be the blow-up at a
point away from $\linf$. $\M _{k}^{n}$ can be alternatively described as the
moduli space of rank $n$ holomorphic bundles on $\til{\cnums \P}^{2}$ with
$c_{1}=0$ and $c_{2}=k$ and with a fixed holomorphic
trivialization on $\linf$.
\end{abstract}

\maketitle

\markboth{Rank Stable Topology of Instantons on $\cpbar $}
{James Bryan and Marc Sanders}
\renewcommand{\sectionmark}[1]{}

\section{Introduction}

In his 1989 paper \cite{Ta}, Taubes studied the stable topology of the
based instanton moduli spaces. He showed that if $\M_{k}^{n}(X)$ denotes
the moduli space of based $SU(n)$-instantons of charge $k$
on $X$, then there is a map $\M _{k}^{n}(X)\to \M _{k+1}^{n}(X)$ and,
in the direct limit topology, $\M _{\infty }^{n}(X)$ has the homotopy type
of $\operatorname{Map}_{0}(X,BSU(n))$.

There is also a map $\M _{k}^{n}(X)\hookrightarrow \M _{k}^{n+1}(X)$ given by
the direct sum of a connection with the trivial connection on a trivial
line bundle and one can consider the direct limit $\M _{k}^{\infty }(X)$.
For the case of
$X=S^{4}$ with the round metric, it was shown by Kirwan and also by Sanders
(\cite{Kir},\cite{Sa}) that the direct limit has the homotopy type of
$BU(k)$.

In this note we consider the case of $X=\cpbar $ where $\cpbar $ denotes
the complex projective plane with the Fubini-Study metric and the opposite
orientation of the one induced by the complex structure. Our result is:
\begin{theorem}\label{thm: main result}
$\M _{k}^{\infty }(\cpbar )$ has the homotopy type of $BU(k)\times BU(k)$.
\end{theorem}

The main tool in the proof of the theorem is a construction of the moduli
spaces $\M _{k} ^{n}(\cpbar )$ due to King \cite{Ki}. In general,
Buchdahl \cite{Bu} has shown that, for appropriate metrics on the
$N$-fold connected sum $\# _{N}\cpbar $, the moduli spaces $\M _{k}^{n}(\#
_{N}\cpbar )$ are diffeomorphic to certain spaces of equivalence classes of
holomorphic bundles on $\cnums \P^{2}$ blown-up at $N$ points. The
universal $U(k)\times U(k) $ bundle that appears giving the homotopy
equivalence of theorem \ref{thm: main result} can be constructed as
higher direct image bundles (see section \ref{section:proof}).

\begin{rem}\label{rem:motivation}
The cofibration $S^2\to \cpbar \to S^4$ gives rise to the fibration
of mapping spaces $\Omega^4BSU(n)\to {\rm Map}_{\ast}(\cpbar , BSU(n))\to
\Omega^2BSU(n)$ which for K-theoretic reasons is a trivial fibration in the
limit over $n$. The total space of this fibration is homotopy equivalent to
the space of based gauge equivalence classes of all connections on
$\cpbar$.  Thus,
from Taubes' result, $\M_k^\infty $
must have
the property that taking the limit over $k$ gives $BU\times BU$.  For $S^4$,
similar remarks
imply that $\lim_{k\to \infty} \M _{k} ^{\infty}(S^4)\simeq BU$ and the
inclusion
of $\M _k ^{\infty} (S^4)$ into this limit has been shown to be (up to
homotopy) the
natural inclusion $BU(k)\hookrightarrow BU$ (\cite {Sa}).
Theorem 1.1
and these results for $S^4$ suggest a general conjecture which is supported
by the fact that the higher direct image bundle giving our homotopy equivalence
generalizes in an appropriate way.
\end{rem}

\begin{conj}\label{conj:N fold connected sum}
For appropriate metrics on $\# _{N}\cpbar $, $\M _{k}^{\infty }(\#
_{N}\cpbar )$ has the homotopy type of a product $BU(k)\times \cdots
\times BU(k)$ with $N+1$ factors.
\end{conj}
\begin{rem}\label{rem:Bott periodicity}
Combining theorem \ref{thm: main result} with Taubes' stabilization result
leads to an alternate proof of Bott periodicity for the unitary group.
There is a natural map from instantons on $S^4$ to those on $\cpbar$ (given by
pull-back) which in the limit as $n\to \infty$ is homotopy equivalent to the
diagonal $BU(k) \to BU(k)\times BU(k)$. Taking the limit as $k\to\infty$ and
applying Taubes' result, the diagonal map appears as the inclusion of fibers
in
the fibration $BU\simeq\Omega_k^4BU\to BU\times BU\to \Omega^2(BSU)$ (see
remark 1.1).
However, the map $BU\rightarrow BU \times BU$ has homotopy fiber $U\times U
\bigg/ U\: \simeq U$,
and therefore, given the above fibration we must have
$U\simeq\Omega^3BSU\simeq\Omega^2 SU\simeq\Omega^2 U$.

Tian (\cite {Ti})
noticed that results for limits of instantons on $S^4$ (see \cite {Kir} or
\cite {Sa})
already imply the four-fold periodicity $\Omega^4BU\simeq Z\times BU$.  The
${\cpbar}$ case thus gives the
finer two-fold periodicity $\Omega^2BU\simeq Z\times BU$, as one might expect
due to the nontrivial $S^2\hookrightarrow {\cpbar}$.

In a future paper we will study limits of $Sp(n)$ and $SO(n)$ instantons
on $\cpbar$ and their relationships to those on $S^4$.  As an amusing
corollary, we will be able to rederive many of the
Bott periodicity
relationships among
$Sp$, $U$, $SO$, and their homogeneous spaces.

\end{rem}

\section{The construction of $\M _{k}^{n}(\cpbar )$}

Let $x_{0}\in \cpbar $ be the base point. Since $\cpbar -\{x_{0}  \}$ is
conformally equivalent to $\til{\cnums} ^{2}$, the complex plane
blown-up at the origin, $\M _{k}^{n}(\cpbar )$ can
be regarded as instantons on $\til{\cnums} ^{2}$ based ``at infinity''.
Buchdahl \cite{Bu} proved an analogue in this non-compact setting of
Donaldson's theorem relating instantons to holomorphic bundles: Let
$\til{\cnums}^{2}_{N}$ be the complex plane blown-up at $N$ points with a
K\"ahler metric.  Then $\til{\cnums}^{2}_{N}$ has a ``conformal
compactification'' to $\#_{N}\cpbar$ and a ``complex compactification'' to
$\til{\cnums \P}^{2}_{N}$ (the projective plane blown-up at $N$ points).
We have added a point
$x_{0}$ in the former case and a complex projective line $\linf$ in the latter.

Define $\msalg{k}{n}(\til{\cnums \P}^{2}_{N})$ to be the moduli space
consisting of pairs $(\E,\tau
)$ where $\E $ is a rank $n$ holomorphic bundle on $\til{\cnums \P}
^{2}_{N}$ with $c_{1}(\E )=0$, $c_{2}(\E )=k$, and where $\tau :\E
|_{\linf}\to \cnums ^{n}\otimes \O _{\linf}$ is a
holomorphic trivialization of $\E $ on $\linf$.

There is a natural map
$\Phi:\M_{k}^{n}(\# _{N}\cpbar )\longrightarrow\msalg{k}{n}(\til{\cnums
\P}^{2}_{N})$ defined as follows. Let
$p:\til{\cnums \P}^{2}_{N}\to \# _{N}\cpbar  $ be the map that collapses
$\linf\mapsto x_{0}$. If
$[A]\in\ms{k}{n}$ then the $\dbar$ operator that defines the holomorphic bundle
$\mathcal{V}=\Phi(A)$ is taken to be $(d_{p^{*}(A)})^{(0,1)}$,
the anti-holomorphic part
of the covariant derivative defined by the pullback of the connection. The
anti-self-duality of $A$ implies that the curvature of $p^{*}(A)$ is a
$(1,1)$-form and so $\dbar ^{2}=0$.

Buchdahl's theorem is then

\begin{theorem} \label{thm:buchdahl}
The map $\Phi:\M_{k}^{n}(\# _{N}\cpbar
)\longrightarrow\msalg{k}{n}(\til{\cnums \P}^{2}_{N})$ is a
diffeomorphism.
\end{theorem}

The case $N=1$ was first proved by King \cite{Ki}. We
now restrict ourselves to that case and simply write
$\ms{k}{n}$ for $\M_{k}^{n}(\cpbar )$ and $\msalg{k}{n}(\til{\cnums \P}^{2})$.

King constructed $\ms{k}{n}$ explicitly in terms of linear algebra data.
We recall his construction. Consider configurations of linear maps:

\begin{picture}(100,75)(-135,0)
\put(20,50){$W_0$}
\put(56,15){\vector(-2,3){21}}
\put(55,5){$V_{\infty}$}
\put(36,52){\vector(1,0){47}}
\put(83,56){\vector(-1,0){47}}
\put(84,50){$W_1$}
\put(86,45){\vector(-2,-3){19}}
\put(35,27){$b$}
\put(56,45){$x$}
\put(48,60){$a_1,a_2$}
\put(81,27){$c$}
\end{picture}

\noindent
where $W_0$, $W_1$ and $V_{\infty }$ are complex vector spaces of
dimensions $k$, $k$, and $n$ respectively.

A configuration $(a_1,a_2,b,c,x)$ is called
{\em integrable} if it satisfies the equation
$$a_1xa_2-a_2xa_1+bc=0.$$

A configuration $(a_1,a_2,b,c,x)$ is {\em non-degenerate} if it satisfies
the following conditions:
$$\forall \hspace{4pt}
(\lambda_1,\lambda_2),(\mu_1,\mu_2)\in\cnums^2 \text{ such that
}\lambda_1\mu_1 + \lambda_2\mu_2=0\text{ and }(\mu_1,\mu_2)\neq(0,0),$$
\begin{eqnarray*}
\notexists v\in W_1\text{ such that } &&\left\{
\begin{array}{lr}
xa_1v=\lambda_1v & (\mu_1 a_1 + \mu_2 a_2)v=0 \\
xa_2v=\lambda_2v & cv=0
\end{array}
\right.  \\
\text{and }\notexists w\in W_0^\star\text{ such that } &&\left\{
\begin{array}{lr}
x^*a_1^*w=\lambda_1w & (\mu_1 a^*_1 + \mu_2 a^*_2)w=0 \\
x^*a_2^*w=\lambda_2w & b^*w=0
\end{array}
\right.
\end{eqnarray*}

Let $A_{k}^{n}$ be the space of all integrable
non-degenerate configurations. $G=Gl(W_{0})\times Gl(W_{1})$ acts
canonically on $A_k^n$. The action is explicitly given by
$$
(g_{0},g_{1})\cdot (a_{1},a_{2},b,c,x)=
(g_{0}a_{1}g_{1}^{-1},g_{0}a_{2}g_{1}^{-1},
g_{0}b,cg_{1}^{-1},g_{1}xg_{0}^{-1})
$$

\begin{theorem}\label{thm: monad description}
The moduli space $\ms{k}{n}$ is isomorphic to
$A_{k}^{n}/G$.
\end{theorem}
\proof
King uses such  configurations
to determine monads that in turn
determine holomorphic bundles. Configurations in the same $G$
orbit determine the same bundle. For the sake of brevity we refer the
reader to \cite{Ki} or \cite{Br} for details. The construction identifies the
vector spaces $W_0$ and $W_{1}$ canonically as $H^1(\mathcal{E}(-\linf))$
and $H^{1}(\mathcal{E}(-\linf+E))$ respectively, where $E\subset \til{\cnums
\P}^{2}$ is the exceptional divisor. The vector space $V_{\infty } $ is
identified with the fiber over $\linf$.

\section{Proof of theorem \ref{thm: main result} }\label{section:proof}

We prove the theorem in two steps: We first show that the space of monad
data $A_{k}^{n}$ forms a principal $G=Gl(k)\times Gl(k)$ bundle over $\M
_{k}^{n}$. We then show that the induced $G$-equivariant inclusion
$A_k^n\hookrightarrow A_k^{n+2k}$ is null-homotopic so that we can conclude
that $A_{k}^{\infty }$ is contractible.

\begin{lem}\label{lem:freeness of G action}
$G$ acts freely on the space of monad data $A_{k}^{n}$.
\end{lem}
\begin{proof}
This is essentially proved in \cite{Ki} where it is implicitly shown that
the non-degeneracy conditions are precisely the conditions that guarantee
freeness. We point out that this also follows more conceptually from the
existence of a universal family $\Bbb{E}\to \ms{k}{n}\times \til{\cnums
\P} ^{2}$ and the cohomological
interpretation of $W_{0}$ and $W_{1}$:

First, the existence of a universal family can be shown via the gauge
theoretic construction:
Let $V$ be a smooth hermitian vector bundle on $\til{\cnums \P}^{2}$ with
$c_{1}(V)=0$ and $c_{2}(V)=k$.  Let $\mathcal{A}^{1,1}_{0}$ denote
unitary connections on $V$
with curvature of pure type $(1,1)$ and that restrict to the trivial
connection on $\linf $ and let $\mathcal{G}_{0}^{\cnums }$ denote the
complex gauge transformations of $V$ that are the identity restricted to
$\linf $.
Then $\ms{k}{n}=\mathcal{A}^{1,1}_{0}/\mathcal{G}^{\cnums }_{0}$. The
quotient
$$(\mathcal{A}^{1,1}_{0}\times V )/\mathcal{G}^{\cnums }_{0}\to \M
_{k}^{n}\times \til{\cnums \P}^{2}$$
will form a universal bundle if the moduli space is smooth and no $\E \in
\M _{k}^{n}$ has non-trivial automorphisms ({\em c.f.} \cite{Fr-Mo} Chapt. IV):

\begin{lem}\label{lem:E has no finite automorphisms, moduli space is smooth}
$\M _{k}^{n}$ is smooth and any $\E \in \M _{k}^{n}$ has no non-trivial
automorphisms preserving $\tau :\E |_{\linf} \to \cnums ^{n}\otimes
\mathcal{O}_{\linf }$.
\end{lem}

By Serre duality $H^{2}(\E \otimes \E^{*} )=H^{0}(\E \otimes \E ^{*}\otimes
K)^{*}$. Since $\E \otimes \E ^{*}$ is trivial on $\linf $, it is
also trivial on nearby lines.  Any section of $\E \otimes \E^{*}
\otimes K$ restricts to a section of $\cnums ^{n^{2}}\otimes
\mathcal{O}_{\linf }(-3)$ and so must
vanish on $\linf $. Likewise, it must vanish on nearby lines and so it is
$0$ on an open set and must be identically $0$. Thus  $H^{2}(\E\otimes \E
^{*})=0$ and smoothness follows once we show there are no automorphisms.

Suppose that there exists an automorphism $\phi \in H^{0}(\E \otimes \E
^{*})$ such that $\phi \neq \1 $ and $\phi $ preserves $\tau $ so that
$\phi |_{\linf }=\1 |_{\linf }$. Then $\phi -\1 $ is a non-zero section of
$\E \otimes \E^{*} $
vanishing on $\linf $. We then get an injection $0\to \O(\linf )\to  \E
\otimes \E ^{*}$. Restricting this sequence to $\linf $ we get an injection
$0\to \O_{\linf }(1)\to \O_{\linf }\otimes \cnums ^{n^{2}}$
which is a contradiction.

Let $\pi :\ms{k}{n}\times \til{\cnums\P }^{2}\to \ms{k}{n}$. The
higher direct image sheaves $R^{1}\pi _{*}(\Bbb{E}(-\linf))$ and $R^{1}\pi
_{*}(\Bbb{E}(-\linf+E))$  are locally free and rank $k$. This follows from
the index theorem and the vanishing of the $H^{0} $ and $H^{2}$
cohomology of $\E (-\linf )$ and $\E (-\linf +E)$. The $H^{0}$ vanishing follows
by again considering the restriction of a section of the bundles to lines
nearby to $\linf $. Using Serre duality and the same argument, one gets the
vanishing for $H^{2}$.

Thus the vector spaces $W_{0}$ and $W_{1}$ are the fibers of the vector
bundles $R^{1}\pi_{*}(\Bbb{E}(-\linf))$ and
$R^{1}\pi_{*}(\Bbb{E}(-\linf+E))$. The $G$-orbit of a configuration giving a
bundle $\E$  can be identified with the group of isomorphisms
$g_0:H^1(\E(-\linf))\to
\cnums^k$ and $g_1:H^{1}(\E(-\linf+E))\to \cnums^k$. Thus
$A_{k}^{n}$ is realized precisely as the total space of the
principal $Gl(k)\times Gl(k)$  bundle associated to
$R^{1}\pi_{*}(\Bbb{E}(-\linf))\oplus R^{1}\pi_{*}(\Bbb{E}(-\linf+E))$.

\end{proof}

Recall that the map $\ms{k}{n}\hookrightarrow \ms{k}{n+1}$ is defined
by the direct sum with the trivial connection: $[A]\mapsto [A\oplus \theta]$.
In terms of holomorphic bundles this is $\E \mapsto \E \oplus \mathcal{O}$.
Tracing through the monad construction, it is easy to see that the
inclusion induces the $G$-equivariant map $A_k^n\hookrightarrow A_k^{n+1}$
given by $(a_{1},a_{2},x,b,c)\mapsto (a_{1},a_{2},x,b',c')$ where
$b^{\prime}$ is $b$ with an extra first column of zeroes and $c^{\prime}$ is
$c$ with an extra first row of zeroes.
Define $A_k^{\infty}$ to be the direct
limit $\lim_{n\to\infty}A_k^n$ so that there is a homeomorphism between $\M
_k^{\infty}$ and $A_k^{\infty}/G$

\begin{lem}
$A_k^{\infty}$ is a contractible space.
\end{lem}

\begin{proof}
Since the $A_k^{n}$'s are
algebraic varieties and the maps
$A_k^n\to A_k^{n+1}$ are algebraic, they admit triangulations compatible with
the maps. Thus $A_k^\infty$ inherits the structure of a CW-complex and so it
is sufficient to show that all of its homotopy groups are zero.
To this end we prove that for any
$k$ and $l$   there is an $r>l$ such that the natural inclusion
from $A_k^n\hookrightarrow A_k^r$ is homotopically trivial.

Consider the homotopy
$H_t:A_k^n\to A_k^{2k+n}$
defined as follows:
\[H_t((a_1,a_2,x,b,c))=
( (1-t)a_1,\; (1-t)a_2,\; (1-t)x,\; b_t,\; c_t)\]
where
\[c_t=
\begin{pmatrix}tI_k\\ 0_{k,k}\\ (1-t)c
\end{pmatrix}
\;\;{\rm ,}\;\;
b_t=(0_{k,k},\; tI_k,\; (1-t)^2b),\]
$I_k$ is the $k\times k$ identity matrix and $0_{k,k}$ is the
$k\times k$ zero matrix.  To see that
$H_t(v)\in A_k^{n+2k}$ for any $v\in A_k^n$, we
check that the  integrability and non-degeneracy conditions
are satisfied for all $0\leq t\leq 1$.
Integrability holds because $b_tc_t=(1-t)^3bc$.
Non-degeneracy is satisfied for all $t\not=0$
because there is a full rank $k\times k$ block, $tI_k$,
in both
$c_t$ and $b_t$.  Furthermore, $H_0$ is just the inclusion
$A_k^n\hookrightarrow A_k^{n+2k}$, so
non-degeneracy also holds when $t=0$.
Finally, note that $H_1$ is a constant map.
\end{proof}

These lemmas show that $A_{k}^{\infty }$ is a contractible space acted on
freely by $G=Gl(k)\times Gl(k)$ and $A_{k}^{\infty }/G=\ms{k}{\infty }$.
Thus $\ms{k}{\infty }$ is homotopic to $BG$ which in turn has the homotopy
type of $BU(k)\times BU(k)$. We end by remarking that the proof shows that the
universal $U(k)\times U(k)$ bundle is the bundle that restricts to any of the
finite $\M_k^n$'s as $R^{1}\pi_{*}(\Bbb{E}(-\linf))\oplus
R^{1}\pi_{*}(\Bbb{E}(-\linf+E))$.

\end{document}